\documentclass[a4paper,11pt]{amsart}

\usepackage{graphicx}

\newtheorem{theorem}{Theorem}[section]
\newtheorem{lemma}[theorem]{Lemma}
\newtheorem{corollary}[theorem]{Corollary}
\newtheorem{proposition}[theorem]{Proposition}
\newtheorem{definition}[theorem]{Definition}

\newtheorem{remark}[theorem]{Remark}

\newcommand{\pnt}{\circle*{5}}
\newcommand{\R}{\mathbb{R}}
\newcommand{\N}{\mathbb{N}}
\newcommand{\A}{\mathcal{A}}
\newcommand{\D}{\mathcal{D}}
\newcommand{\I}{I}%{\mathcal{I}}
\newcommand{\IH}{\hat{\I}}
\newcommand{\OS}[1]{\mathcal{OS}(#1)}
\DeclareMathOperator{\spn}{span}
\DeclareMathOperator{\rank}{rank}
\DeclareMathOperator{\cl}{cl}
\DeclareMathOperator{\sym}{Sym}
\DeclareMathOperator{\path}{path}
\DeclareMathOperator{\supp}{supp}
\DeclareMathOperator{\hgt}{ht}

\begin{document}
\author{Axel Hultman}
\address{Department of Mathematics, Link\"oping University, SE-581 83, Link\"oping, Sweden.}
\title[Root ideal arrangements]{Supersolvability and the Koszul property of root ideal arrangements}

\begin{abstract}
A {\em root ideal arrangement} $\A_\I$ is the set of reflecting hyperplanes
corresponding to the roots in an order ideal $\I\subseteq \Phi^+$ of the root
poset on the positive roots of a finite crystallographic root system
$\Phi$. A characterisation of supersolvable root ideal arrangements is
obtained. Namely, $\A_\I$ is supersolvable if and only if $\I$ is {\em
  chain peelable}, meaning that it is possible to reach the empty
poset from $\I$ by in each step removing a maximal chain which is also
an order filter. In particular, supersolvability is preserved under
taking subideals. We identify the minimal ideals that correspond to
non-supersolvable arrangements. There are essentially two
such ideals, one in type $D_4$ and one in type $F_4$. By showing that
$\A_\I$ is not line-closed if $\I$ contains one of these, we deduce that the Orlik-Solomon algebra
$\OS{\A_\I}$ has the Koszul property if and only if $\A_\I$ is
supersolvable.
\end{abstract}

\maketitle

\section{Introduction}
Let $\Phi$ be a finite, crystallographic root system with set
of positive roots $\Phi^+$. The orthogonal complements of the elements
of $\Phi^+$ form a hyperplane arrangement known as the {\em Weyl
  arrangement}, or {\em Coxeter arrangement}, of $\Phi$; this is the set of reflecting hyperplanes
of the corresponding reflection group.

With any set of positive roots $R\subseteq \Phi^+$ is associated the
corresponding subarrangement $\A_R$ of the Weyl arrangement. That is,
\[
\A_R = \{H_\gamma\mid \gamma\in R\},
\]
where $H_\gamma$ is the orthogonal complement of $\gamma$.

In \cite{ST}, Sommers and Tymoczko initiated the study of an interesting
family of subarrangements of the Weyl arrangement. Suppose $\I
\subseteq \Phi^+$ is an order ideal in the root poset on
$\Phi^+$. That is, if $\gamma_1\in \I$ and $\gamma_2 < \gamma_1$,
then $\gamma_2\in \I$. It was shown in
\cite{ST} for $\Phi$ of classical type (and conjectured in general)
that every {\em root ideal arrangement} $\A_\I$ is free in the sense of Terao
\cite{terao}. This result was extended to arbitrary finite
crystallographic type by Abe {\em et al.\ }\cite{ABCHT}. 

%Being free, a root ideal arrangement
%has associated {\em exponents} that show up when factoring its
%characteristic polynomial. These exponents also govern (party
%conjecturally) a factorisation of the Poincar\'e polynomial of a Hessenberg variety
%indexed by $\I$. This 
%generalises how the ordinary Weyl group exponents appear both in the characteristic
%polynomial of the Weyl arrangement and the Poincar\'e
%polynomial of the corresponding flag variety; this situation is
%recovered when $\I = \Phi^+$. 

A stronger property than freeness is supersolvability. Although being
free, root ideal arrangements are not in general
supersolvable. Indeed, for irreducible root systems, it is known
that $\A_{\Phi^+}$ is supersolvable if and only if $\Phi$ is of type $A_n$,
$B_n$, $C_n$ or $G_2$ \cite{BI}.

Any matroid $\A$ has an associated Orlik-Solomon algebra $\OS{\A}$, as
introduced in \cite{OS}. We refer to \cite{OT} or \cite{yuzvinsky} for
accounts of the many attractive properties of $\OS{\A}$. Let us here
merely recall that if $\A$ is a complex hyperplane arrangement (meaning that a flat of
the matroid is the collection of all hyperplanes containing a given
subspace), $\OS{\A}$ is isomorphic to the cohomology algebra of the complement of $\cup\A$. 

Shelton and Yuzvinsky \cite{SY} proved that $\OS{\A}$ is a Koszul algebra
whenever $\A$ is supersolvable. Whether the converse holds remains an
open question which has been answered affirmatively for certain classes of
arrangements. For example, it is known for hypersolvable
arrangements \cite{JP}, graphic arrangements \cite{SS} and arrangements whose minimal broken circuits
are pairwise disjoint \cite{LR}.

Let us call a poset chain peelable if it is possible to remove
all its elements by using moves that remove, in each step, a set of
elements that simultaneously form a maximal
totally ordered subset and an order filter (i.e.\ the
complement of an order ideal). The first main result of the present
paper is the following characterisation of 
supersolvable root ideal arrangements:
\begin{theorem} \label{th:main_peelable}
Suppose $\Phi$ is a finite crystallographic root system and that
$\I\subseteq \Phi^+$ is an order ideal in the root poset on the
positive roots. Then, the root ideal arrangement $\A_\I$ is
supersolvable if and only if the poset $\I$ is chain peelable.
\end{theorem}

It is clear that chain peelability of a poset is inherited by all its
ideals. This means that there exist minimal ``bad'' root poset ideals in the
sense that $\A_\I$ is not supersolvable if and only if $\I$ contains one
of the bad ideals as a subset. We identify these bad ideals, finding
that there are essentially only two of them. One is generated by the roots of
height $3$ in type $D_4$, whereas the other is generated by the roots of
height $4$ in type $F_4$.

Falk \cite{falk} introduced the concept of line-closedness and showed
for a matroid $\A$ that having this property is a
necessary condition for $\OS{\A}$ being quadratic. This, in turn, is a
necessary condition for $\OS{\A}$ being Koszul; see e.g.\ \cite{BGS}. By showing that a
root ideal arrangement $\A_\I$ cannot be line-closed if $\I$ contains a
bad ideal, we deduce
\begin{theorem}\label{th:main_koszul}
A root ideal arrangement is supersolvable if and only if its
Orlik-Solomon algebra is Koszul.
\end{theorem}
Thus, root ideal arrangements form yet another class of
arrangements for which the Koszul property is equivalent to
supersolvability.

The remainder of this paper is organised as follows. In the next
section, we review some properties of root posets and agree on root system
notation. In Section \ref{se:arrangements}, conventions on
arrangements are established and a convenient characterisation of
supersolvability is recalled. We then start pursuing Theorem
\ref{th:main_peelable}. General restrictions on the structure of
supersolvable root ideal arrangements 
are collected in Section \ref{se:structure}. Chain
peelable posets are formally defined in
Section \ref{se:peelability} where it is shown that that $\A_\I$ is
supersolvable if $\I$ is chain peelable. The converse is
established for simply laced $\Phi$ in Section \ref{se:simply_laced},
whereas multiply laced root systems are handled in Section
\ref{se:multiply_laced}. In the final section, the equivalence between
the Koszul property of $\OS{\A_\I}$ and supersolvability of $\A_\I$ is
deduced. 

\section{Root systems and root posets}

Consider a finite, crystallographic root system $\Phi$ with set of
positive roots $\Phi^+ \subset \Phi$. Denote by $\Delta \subseteq \Phi^+$ the set
of simple roots. Then, $\Delta$ is a basis for the ambient Euclidean
space $V \cong \R^{\rank \Phi}$. For $\gamma\in \Phi^+$ and $\alpha\in
\Delta$, let $\gamma_\alpha$ be the corresponding coordinate. Thus,
\[
\gamma = \sum_{\alpha \in \Delta} \gamma_\alpha \alpha, 
\]
where all $\gamma_\alpha$ are nonnegative integers. The {\em height} of $\gamma$ is
defined by 
\[
\hgt(\gamma) = \sum_{\alpha\in \Delta}\gamma_\alpha
\]
and the {\em support} of $\gamma$ is 
\[
\supp\gamma = \{\alpha\in \Delta\mid \gamma_\alpha \neq 0\}.
\]
This support induces a connected subgraph of the Dynkin
diagram of $\Delta$ and, conversely, the sum of the nodes of any
connected subgraph of the Dynkin diagram is a positive root.

The reflections through the orthogonal complements of the roots are
the reflections of a finite reflection group. Let $s_\alpha$ denote the reflection
corresponding to the simple root $\alpha\in \Delta$. Then, $s_\alpha$
permutes $\Phi^+\setminus \{\alpha\}$ whereas
$s_\alpha(\alpha) = -\alpha$.

The {\em root poset} on $\Phi^+$ is the partial order defined by
$\gamma \le \gamma'$ if and only if $\gamma_\alpha \le
\gamma'_\alpha$ for all $\alpha \in \Delta$. The cover relation $\lhd$ 
satisfies $\hgt(\gamma) = \hgt(\gamma') - 1$ whenever $\gamma \lhd
\gamma'$. In other words, $\gamma \lhd \gamma'$ means
$\gamma'-\gamma \in \Delta$. The root posets of type $D_4$ and $F_4$
are depicted in Figures \ref{fi:D4} and \ref{fi:F4}, respectively.

For any $S \subseteq \Phi$, $\Phi\cap \spn S$ is a root subsystem of $\Phi$. We
collect here convenient properties of such subsystems that shall
be used frequently, often without explicit mentioning. The first
observation follows, for example, by considering root lengths and
using the fact that every root is in the orbit of a simple root under
the action of the associated reflection group.
\begin{lemma}
If $\Phi$ has a doubly (respectively, triply) laced root subsystem,
then $\Phi$ is doubly (respectively, triply) laced.
\end{lemma}
In particular, every root subsystem of a simply laced system is itself
simply laced.

\begin{lemma}\label{le:ranktwo}
Any rank two subsystem of $\Phi$ contains at most one pair of positive roots
that are incomparable in the root poset on $\Phi^+$. If they exist, these two
elements are minimal among all positive roots in the subsystem.
\begin{proof}
Suppose $\Psi \subseteq \Phi$ is a rank two subsystem. Let
$\gamma_1,\gamma_2$ denote the simple roots of $\Psi$ that correspond to the
choice of positive roots $\Psi^+ = \Phi^+
\cap \Psi$. The other
elements of $\Psi^+$ are positive linear combinations of the
$\gamma_i$, hence larger than $\gamma_i$. Inspection of all rank two
root systems reveals that $\Psi^+ \setminus\{\gamma_1,\gamma_2\}$ is
totally ordered. The assertion follows.
\end{proof}
\end{lemma}
Note that $\gamma_1$ and $\gamma_2$ in the above proof may be comparable
in the root poset on $\Phi^+$. Hence, it is not in general possible to
replace ``at most'' by ``exactly'' in the statement of Lemma
\ref{le:ranktwo}.

Let $(\cdot, \cdot)$ denote the Euclidean inner product on $V$. By
inspecting all simply and doubly laced rank two root systems (namely, 
$A_1\times A_1$, $A_2$ and $B_2 = C_2$), one readily verifies the
following lemma.
\begin{lemma}\label{le:rank2}
Suppose $\Phi$ is a finite crystallographic root system without triple
bonds in the Dynkin diagram. For any two distinct positive roots $\beta, \gamma\in \Phi^+$,
the following assertions hold:
\begin{itemize}
\item If $(\beta, \gamma) > 0$, then $\beta - \gamma \in \Phi$ and
  $\beta + \gamma \not \in \Phi$.
\item If $(\beta, \gamma) < 0$, then $\beta - \gamma \not \in \Phi$ and
  $\beta + \gamma \in \Phi$.
\item If $(\beta, \gamma) = 0$ and $\Phi$ is simply laced, then $\beta
  - \gamma \not \in \Phi$ and $\beta + \gamma \not \in \Phi$.
\end{itemize}
\end{lemma}

\begin{figure}[t]
\scalebox{.8}{
\begin{picture}(320,220)(-10,-5)
\put (0,0){\pnt}\put (10,-5){$1000$}   
\put (100,0){\pnt}\put (110,-5){$0100$}   
\put (200,0){\pnt}\put (210,-5){$0010$}   
\put (300,0){\pnt}\put (310,-5){$0001$}
\put (50,50){\pnt}\put (60,47){$1100$}   
\put (150,50){\pnt}\put (160,47){$0110$}   
\put (250,50){\pnt}\put (260,47){$0101$}
\put (50,100){\pnt}\put (60,97){$1110$}
\put (150,100){\pnt}\put (160,97){$1101$}
\put (250,100){\pnt}\put (260,97){$0111$}
\put (150,150){\pnt}\put (160,150){$1111$}
\put (150,200){\pnt}\put (160,197){$1211$}

\put (0,0){\line(1,1){50}}
\put (100,0){\line(1,1){50}}
\put (100,0){\line(3,1){150}}
\put (100,0){\line(-1,1){50}}
\put (200,0){\line(-1,1){50}}
\put (300,0){\line(-1,1){50}}

\put (50,50){\line(0,1){50}}
\put (50,50){\line(2,1){100}}
\put (250,50){\line(-2,1){100}}
\put (250,50){\line(0,1){50}}
\put (150,50){\line(-2,1){100}}
\put (150,50){\line(2,1){100}}

\put (150,100){\line(0,1){50}}
\put (50,100){\line(2,1){100}}
\put (250,100){\line(-2,1){100}}

\put (150,150){\line(0,1){50}}

\put (0,170){\circle*{8}}
\put (30,170){\circle*{8}}
\put (30,200){\circle*{8}}
\put (60,170){\circle*{8}}
\put (0,170){\line(1,0){30}}
\put (30,170){\line(1,0){30}}
\put (30,170){\line(0,1){30}}

\put (-5, 158){$\alpha_1$}
\put (25, 158){$\alpha_2$}
\put (55, 158){$\alpha_3$}
\put (37, 197){$\alpha_4$}
\end{picture}
}
\caption{The type $D_4$ Dynkin diagram and root poset. The elements of
the root poset are labelled by coordinate sequences. For example, the
label ``$0101$'' indicates the root $\alpha_2 + \alpha_4$.} \label{fi:D4}
\end{figure}
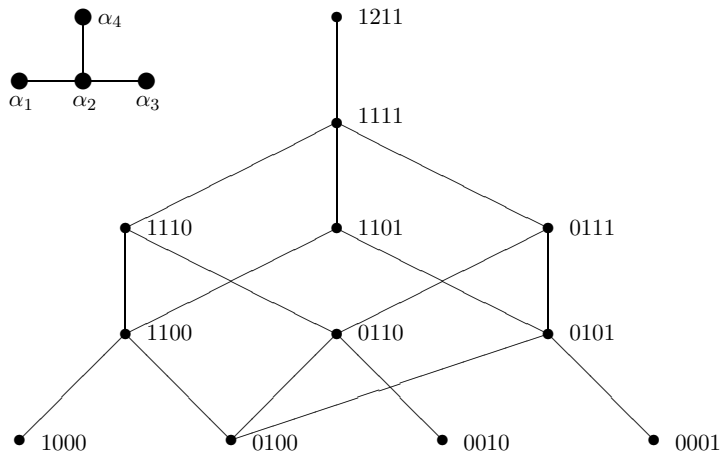

\begin{figure}[t]
\scalebox{.65}{
\begin{picture}(360,540)(-30,-25)
\put (0,0){\pnt}\put (10,-3){$1000$}   
\put (100,0){\pnt}\put (110,-3){$0100$}   
\put (200,0){\pnt}\put (210,-3){$0010$}   
\put (300,0){\pnt}\put (310,-3){$0001$}   
\put (50,50){\pnt}\put (60,47){$1100$}   
\put (150,50){\pnt}\put (160,47){$0110$}   
\put (250,50){\pnt}\put (260,47){$0011$}   
\put (50,100){\pnt}\put (60,97){$1110$}
\put (150,100){\pnt}\put (160,97){$0210$}
\put (250,100){\pnt}\put (260,97){$0111$}
\put (50,150){\pnt}\put (60,147){$1210$}
\put (150,150){\pnt}\put (160,150){$1111$}
\put (250,150){\pnt}\put (260,150){$0211$}
\put (50,200){\pnt}\put (60,197){$2210$}
\put (150,200){\pnt}\put (160,197){$1211$}
\put (250,200){\pnt}\put (260,197){$0221$}
\put (100,250){\pnt}\put (110,247){$2211$}
\put (200,250){\pnt}\put (210,247){$1221$}
\put (100,300){\pnt}\put (110,297){$2221$}
\put (200,300){\pnt}\put (210,297){$1321$}
\put (150,350){\pnt}\put (160,347){$2321$}
\put (150,400){\pnt}\put (160,397){$2421$}
\put (150,450){\pnt}\put (160,447){$2431$}
\put (150,500){\pnt}\put (160,497){$2432$}
\put (0,0){\line(1,1){50}}
\put (100,0){\line(1,1){50}}
\put (200,0){\line(1,1){50}}
\put (100,0){\line(-1,1){50}}
\put (200,0){\line(-1,1){50}}
\put (300,0){\line(-1,1){50}}
\put (50,50){\line(0,1){50}}
\put (150,50){\line(0,1){50}}
\put (250,50){\line(0,1){50}}
\put (150,50){\line(-2,1){100}}
\put (150,50){\line(2,1){100}}
\put (150,100){\line(-2,1){100}}
\put (150,100){\line(2,1){100}}
\put (50,100){\line(0,1){50}}
\put (250,100){\line(0,1){50}}
\put (50,100){\line(2,1){100}}
\put (250,100){\line(-2,1){100}}
\put (50,150){\line(0,1){50}}
\put (150,150){\line(0,1){50}}
\put (250,150){\line(0,1){50}}
\put (50,150){\line(2,1){100}}
\put (250,150){\line(-2,1){100}}
\put (50,200){\line(1,1){50}}
\put (150,200){\line(1,1){50}}
\put (150,200){\line(-1,1){50}}
\put (250,200){\line(-1,1){50}}
\put (100,250){\line(0,1){50}}
\put (200,250){\line(0,1){50}}
\put (200,250){\line(-2,1){100}}
\put (100,300){\line(1,1){50}}
\put (200,300){\line(-1,1){50}}
\put (150,350){\line(0,1){50}}
\put (150,400){\line(0,1){50}}
\put (150,450){\line(0,1){50}}

\put (-30,400){\circle*{8}}
\put (0,400){\circle*{8}}
\put (30,400){\circle*{8}}
\put (60,400){\circle*{8}}
\put (-30,400){\line(1,0){30}}
\put (30,400){\line(1,0){30}}
\put (00,402){\line(1,0){30}}
\put (00,398){\line(1,0){30}}
\put (9,400){\line(2,-1){10}}
\put (9,400){\line(2,1){10}}
\put (-35, 387){$\alpha_1$}
\put (-5, 387){$\alpha_2$}
\put (25, 387){$\alpha_3$}
\put (55, 387){$\alpha_4$}
\end{picture}
}
\caption{The type $F_4$ Dynkin diagram and root poset. The root labels
are coordinate sequences, as in Figure \ref{fi:D4}.} \label{fi:F4}
\end{figure}
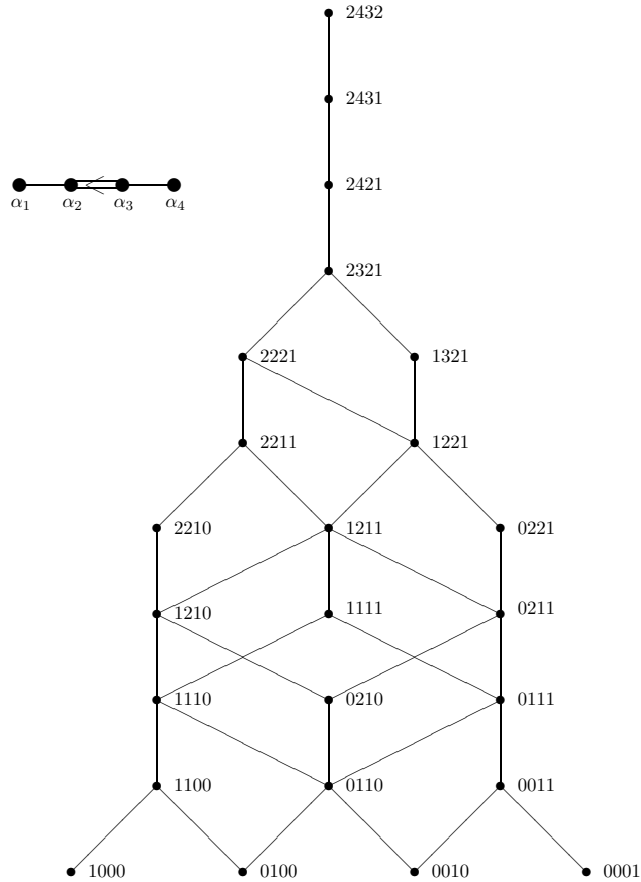

\section{Arrangements and supersolvability} \label{se:arrangements}

Like any vector collection, a subset $R \subseteq \Phi^+$ determines a matroid whose bases are the
maximal linearly independent subsets of $R$. From a matroidal point
of view, $R$ is equivalent to the hyperplane arrangement $\A_R$ consisting of the
orthogonal complements of the roots in $R$. In what follows, we shall
only consider properties of hyperplane arrangements that are
completely determined by the matroidal structure. Therefore we can,
and shall, identify any root ideal arrangement $\A_\I$ with its root
ideal $\I \subseteq \Phi^+$. In this spirit, an {\em
  arrangement} from now on is a collection $\A$ of pairwise non-parallel vectors in some Euclidean space
$V \cong \R^n$. For convenience, we shall assume that $\A$ is {\em essential}, i.e.\ that it
spans $V$.

Let $\A$ be an arrangement. A {\em flat} of $\A$ is a closed set with
respect to the closure operator defined by $\cl(S) = \A \cap \spn S$ for $S\subseteq \A$. The
{\em rank} of a flat (or any other subset of $\A$) is the dimension of
its span.

An arrangement is {\em supersolvable} if its lattice of flats
is supersolvable in the sense of Stanley
\cite{stanley}, i.e.\ if there is a maximal chain of modular 
flats. Bj\"orner, Edelman and Ziegler \cite[Theorem 4.3]{BEZ} showed
that the following alternative definition, which is better suited for
our purposes, is equivalent. For brevity, we refer to a flat of rank
$k$ as a {\em $k$-flat}.

\begin{definition}\label{de:supersolvable}
An arrangement $\A$ is {\em supersolvable} if it can be partitioned
$\A = \Pi_1\biguplus \cdots \biguplus \Pi_n$ in such a way that
for all $i\in \{1, \ldots, n\}$ the following conditions hold:
\begin{itemize}
\item The rank of the subarrangement $\A_i = \Pi_1\biguplus \cdots \biguplus
\Pi_i$ is $i$.
\item No $2$-flat of $\A_i$ is a subset of $\Pi_i$.
\end{itemize}
\end{definition}

Let us say that the expression $\Pi_1\biguplus \cdots \biguplus \Pi_n$ is a
{\em supersolving partition} of $\A$ if it satisfies the conditions of
Definition \ref{de:supersolvable}.

If $\A$ is supersolvable, $\A$ is {\em free} in the sense of Terao
\cite{terao}, meaning that the $\sym V^*$-module $\D(\A)$ of
derivations is free. In this case, the multiset of cardinalities of
the blocks of any supersolving partition coincides with the multiset of
degrees of any homogeneous basis for $\D(\A)$; these degrees are known
as the {\em exponents} of $\A$. Recall from the introduction that all
root ideal arrangements are free. The corresponding exponents seem to
play a role for the Poincar\'e polynomials of Hessenberg varieties
analogous to (and generalising) that which is played by the ordinary
exponents of Weyl groups for the corresponding flag varieties; see
\cite{ST} for more details.

\section{The structure of supersolving partitions}\label{se:structure}

Let $\I$ be a root poset ideal. The next lemma shows that any block of a supersolving partition of
$\I$ must be of one of two specific kinds of sets that we now
define. The first kind is indexed by a simple root $\alpha\in
\Delta$. Let 
\[
F_\I(\alpha) = \{\gamma \in \I \mid \gamma \ge \alpha\}
\]
be the order filter in $\I$ generated by $\alpha$. For the second
kind, choose simple roots $\alpha, \beta\in \Delta$ and
positive integers $a,b$. Define 
\[
G_\I(\alpha, \beta, a, b) = \{\gamma\in \I \mid \gamma_\alpha\alpha +
\gamma_\beta \beta \neq k(a\alpha + b\beta) \text{ for all }k\in \N\}.
\]
\begin{lemma} \label{le:twocases}
Suppose that $\I = \Pi_1\biguplus \cdots \biguplus
\Pi_n$ is a supersolving partition. Then either (a) there is some $\alpha\in
\Delta$ such that $\Pi_n = F_\I(\alpha)$, or (b) there are $\alpha, \beta\in
\Delta$ and positive integers $a,b$ such that $\Pi_n = G_\I(\alpha,
\beta, a, b)$. 

Moreover, if (a) holds, $F_\I(\alpha)$ is a chain. If
(b) holds, $a\alpha + b\beta \in \Phi^+$.
\begin{proof}
Define $d = |\Pi_n\cap \Delta|$ and $\I_{n-1} = \I \setminus
\Pi_n$. Since $\I_{n-1}$ has rank
$n-1$, we have $d \ge 1$. It is not possible to have $d \ge 3$, because
in that case $\Pi_n$ would contain two orthogonal simple
roots and, hence, a $2$-flat. Thus, $d=1$ or $d=2$. 

First assume $d=1$; say $\Pi_n\cap \Delta = \{\alpha\}$. Thus,
$\I_{n-1}\subseteq \spn (\Delta\setminus\{\alpha\})$, because 
the dimension of this span is $n-1$. Therefore, $F_\I(\alpha)
\subseteq \Pi_n$. If equality does not hold, we find $\gamma \in \Pi_n
\setminus F_\I(\alpha)$. Since $\gamma$ and $\alpha$ are incomparable,
all roots in $\I \cap \spn\{\alpha, \gamma\}$ except $\gamma$ belong
to $F_\I(\alpha)$, but then $\Pi_n$ contains an entire
$2$-flat, a contradiction. We conclude that $\Pi_n =
F_\I(\alpha)$ so that situation (a) is at hand. Being an order filter, $F_\I(\alpha)$ cannot contain incomparable
elements since, by Lemma \ref{le:ranktwo}, two such elements would
generate a $2$-flat entirely contained in $F_\I(\alpha)$. Hence,
$F_\I(\alpha)$ is in fact a chain in this case. 

Now consider the case $d=2$ with $\Pi_n \cap \Delta =
\{\alpha,\beta\}$. No $2$-flat is contained in $\Pi_n$, so $a\alpha +
b\beta\in \I_{n-1}$ for some positive integers $a$ and $b$. (In
particular, $a \alpha + b \beta$ is a root.) The span
of $\{a\alpha + b\beta\}\cup \Delta \setminus \{\alpha, \beta\}$ has
dimension $n-1$, hence must be equal to $\spn \I_{n-1}$. Thus,
$\gamma_\alpha\alpha + \gamma_\beta\beta$ must be parallel to $a\alpha
+ b\beta$ for every $\gamma \in \I_{n-1}$. This shows $\Pi_n \supseteq G_\I(\alpha,
\beta, a, b)$. Equality must in fact hold, for if $\gamma_\alpha\alpha +
\gamma_\beta\beta$ is parallel to $a\alpha + b\beta$, then $\gamma$ is the
only positive root with this property in the $2$-flat generated by $\alpha$ and
$\gamma$, forcing $\gamma\in \I_{n-1}$. Thus (b) is satisfied.
\end{proof}
\end{lemma}

\begin{remark}
There are just a few possible values that $a$ and $b$ can assume in case
(b) of Lemma \ref{le:twocases}. Since $a\alpha+b\beta \in \Phi^+$,
there is a bond between $\alpha$ 
and $\beta$ in the Dynkin diagram. If the bond is simple, the only
possibility is $(a,b)=(1,1)$. If the bond is double, $\alpha$ and
$\beta$ have different lengths. If $\alpha$ is the long root, either
$(a,b)=(1,1)$ or $(a,b) = (1,2)$.
\end{remark}

Since by construction $F_\I(\alpha)$ is an order filter, $\I\setminus
F_\I(\alpha)$ is an ideal in the root poset. Notice however, that $\I
\setminus G_\I(\alpha, \beta, a, b)$ is not in general an ideal. Thus,
Lemma \ref{le:twocases} seems to suggest that in order to
study supersolvable root ideals one is forced to leave the realm of
ideals. Strictly speaking, this is true; there are root ideals with
supersolvable partitions whose initial segments are not all
ideals.\footnote{For example, if $\I = \Phi^+$ in type $A_2$ with simple roots
$\alpha_1$ and $\alpha_2$, we may choose
$\Pi_2 = \{\alpha_1, \alpha_2\} = G_\I(\alpha_1,\alpha_2,1,1)$. Then,
$\Pi_1 = \{\alpha_1 + \alpha_2\}$ is not an
order ideal even though $\Pi_1\biguplus \Pi_2$ is
supersolving.} The following lemma shows that these segments {\em are} however
ideals in root subsystems.

\begin{lemma} \label{le:ideal}
If $\alpha, \beta\in \Delta$ and $a\alpha+b\beta\in \Phi^+$ for
positive $a$ and $b$, then $\Phi'^+ = \Phi^+ \setminus
G_{\Phi^+}(\alpha,\beta,a,b)$ is the set of positive roots in a root
subsystem $\Phi'$ of $\Phi$ with set of simple roots $\Delta' = \{a\alpha + b\beta\}\cup
\Delta \setminus \{\alpha, \beta\}$. Moreover, the root poset on
$\Phi'^+$ coincides with the induced subposet of the root poset of
$\Phi^+$. In particular, $\I \setminus G_\I(\alpha,\beta,a,b)$ is an
ideal in the root poset on $\Phi'^+$.
\begin{proof}
Define $\Phi' = \Phi \cap \spn \Delta'$. Being the intersection of
$\Phi$ with a hyperplane spanned by roots, $\Phi'$ is a root
subsystem of rank one less than $\Phi$. The elements of $\Phi'^+ = \Phi^+\cap \Phi'$ are positive linear
combinations of the elements of $\Delta'$. Hence, $\Delta'$ can be
chosen as the simple root set of $\Phi'$. With this choice, the order
relation of the root poset on $\Phi'^+$ is precisely that induced from $\Phi^+$.
\end{proof}
\end{lemma}

\section{Chain peelings and supersolvability}\label{se:peelability}

In this section we define chain peelings of posets. First, we
establish a key lemma about the structure of totally ordered
intervals in root posets. From it, we deduce that chain peelings of
$\I$ are supersolving partitions.

For positive roots $\beta_1, \beta_2$, employ the following
notation for intervals in the root poset:
\[
[\beta_1,\beta_2] = \{\gamma\in \Phi^+\mid \beta_1\le \gamma \le \beta_2\}.
\]

A {\em chain} in a poset is a totally ordered subposet.

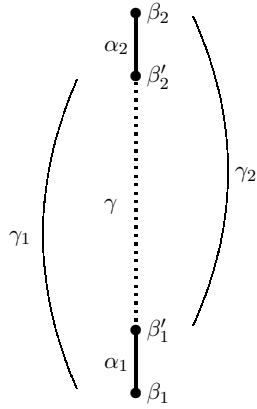
\begin{figure}[t]
\scalebox{0.8}{
\begin{picture}(120,200)(-20,0)
\linethickness{1pt}
\put (50,10){\pnt}\put (55,7){$\beta_1$}   
\put (50,40){\pnt}\put (55,37){$\beta_1'$}   
\put (50,160){\pnt}\put (55,157){$\beta_2'$}   
\put (50,190){\pnt}\put (55,187){$\beta_2$}   

\put (50,10){\line(0,1){30}} \put (35,22){$\alpha_1$}   
\put (50,160){\line(0,1){30}} \put (35,172){$\alpha_2$}   

\multiput (50,44)(0,4){29}{\line(0,1){1}} \put (35,97){$\gamma$}   

\put (97,112){$\gamma_2$}   
\put (-10,82){$\gamma_1$}

\linethickness{.25pt}
\qbezier(77,42)(110,115)(77,188)
\qbezier(22,12)(-10,85)(22,158)  
\end{picture}
}
\caption{Illustration for Lemma \ref{le:chainroot}.} \label{fi:chain}
\end{figure}

\begin{lemma}\label{le:chainroot}
Suppose $[\beta_1, \beta_2]$ is a nonempty chain in the root poset
for distinct positive roots $\beta_1,\beta_2 \in \Phi^+$. Then, $\beta_2 - \beta_1 = k\beta$
for some $\beta \in \Phi^+$ and $k\in \{1,2,3\}$. If $k=3$, the Dynkin
diagram has a triple bond. If $k=2$, the diagram has a double or a
triple bond.
\begin{proof}
By inspecting all rank two root systems, one readily checks that
$\beta_2 - \beta_1 \in k\Phi^+$ can happen only if the Dynkin diagram of the rank two root system
$\Phi^+ \cap \spn \{\beta_1, \beta_2\}$ has an $\ell$-fold bond for
some $\ell \ge k$. This observation verifies the concluding two sentences
of the lemma.

Assume without loss of generality that $\Phi$ is irreducible. In type
$G_2$, the assertion is readily checked, so we may furthermore assume
that there are no triple bonds. Thus, Lemma \ref{le:rank2}
is at our disposal. We shall use it frequently throughout the proof
without explicit mentioning.

If $\beta_2$ covers $\beta_1$, $\beta_2-\beta_1$ is a simple root. If
not, we have $\beta_1 \lhd \beta_1'\leq \beta_2' \lhd \beta_2$. Let
$\alpha_1 = \beta_1'-\beta_1$, $\alpha_2 =
\beta_2-\beta_2'$, $\gamma = \beta_2'-\beta_1'$, $\gamma_1 = \gamma+
\alpha_1$ and, finally, $\gamma_2 = \gamma + \alpha_2$. See Figure
\ref{fi:chain} for an illustration.

Proceed by induction on the length of the chain $[\beta_1,
\beta_2]$. The induction hypothesis shows $\gamma_1,
\gamma_2\in \Phi^+\cup 2\Phi^+$ and $\gamma \in  \Phi^+\cup
2\Phi^+\cup \{0\}$. Moreover, $\alpha_1,\alpha_2\in \Delta$.

There are two cases to consider, depending on whether $\alpha_1$ and
$\alpha_2$ coincide. First, suppose $\alpha_1\neq \alpha_2$. Since $[\beta_1,
\beta_2]$ is a chain, $\beta_2-\alpha_1$ is not a root. Thus,
$(\beta_2, \alpha_1)\le 0$. Similarly, $(\beta_1, \alpha_2)\ge 0$.

Since $\gamma_2 = \gamma_1 + \alpha_2 - \alpha_1$ and the $\alpha_i$
are distinct, it is not possible to have $\gamma_1,\gamma_2\in
2\Phi^+$. Let us assume $\gamma_1\in \Phi^+$; the assumption
$\gamma_2\in \Phi^+$ admits a completely analogous proof.

If $(\gamma_1,\alpha_2) < 0$, $\alpha_2+\gamma_1\in \Phi^+$ and we
are done. Otherwise, the fact that $\beta_2-\gamma_1 = \beta_1+\alpha_2\not\in
\Phi^+$ yields
\[
0 \ge (\gamma_1,\beta_2) = (\gamma_1,\beta_1)+|\gamma_1|^2 +
(\gamma_1,\alpha_2) \ge (\gamma_1,\beta_1)+|\gamma_1|^2.
\]
This is only possible if $\gamma_1$ is short and $\beta_1$ is long. Therefore,
\[
(\beta_1,\beta_2) = |\beta_1|^2 + (\beta_1, \gamma_1) +
(\beta_1,\alpha_2) \ge |\beta_1|^2 + (\beta_1, \gamma_1) > 0,
\]
which leads to $\beta_2-\beta_1 \in \Phi^+$, concluding the proof when
$\alpha_1\neq \alpha_2$.

Now suppose $\alpha_1 = \alpha_2 = \alpha$. If $\gamma = 0$,
$\beta_2-\beta_1 = 2\alpha$ and we are
done. Otherwise, $\supp \gamma\setminus \{\alpha\} \neq \emptyset$.

We have $(\beta_1,\alpha)\le 0$ and $(\beta_2,\alpha)\ge 0$
since $\beta_1+\alpha \in \Phi^+$ and $\beta_2-\alpha\in
\Phi^+$. Thus, $s_\alpha(\beta_1) = \beta_1 + k_1\alpha$ and
$s_\alpha(\beta_2) = \beta_2 - k_2\alpha$ for some $k_1,k_2\in
\{0,1,2\}$. It follows that 
\[
\beta_1\le s_\alpha(\beta_1) < s_\alpha(\beta_2) \le \beta_2.
\]
If one of the weak inequalities is strict,
$s_\alpha(\beta_2-\beta_1) \in (\Phi^+\setminus\{\alpha\})\cup 2
(\Phi^+\setminus\{\alpha\})$ by the induction hypothesis. Then,
$\beta_2-\beta_1  \in (\Phi^+\setminus\{\alpha\})\cup 2 
(\Phi^+\setminus\{\alpha\})$ and we are done. If, on the other hand,
neither of the weak inequalities is strict, we have
\[
s_\alpha(\gamma) = s_\alpha(\beta_2-\beta_1-2\alpha) =
\beta_2-\beta_1+ 2\alpha = \gamma + 4\alpha.
\]
Hence $\gamma\not \in \Phi^+$. By the induction hypothesis, this
implies 
$\gamma = 2\gamma', \gamma'\in \Phi^+$, and $s_\alpha(\gamma') =
\gamma'+2\alpha\in \Phi^+$. Thus, $\gamma'+\alpha\in
\Phi^+$. Observing $\beta_2-\beta_1 = 2(\gamma'+\alpha)$ concludes the proof.
\end{proof}
\end{lemma}

\begin{definition}\label{de:peeling}
Suppose $P$ is a finite poset. If $P$ is a chain, $P$ is a {\em chain peeling} of
itself. Otherwise a {\em chain peeling} of $P$ is a partition of its
elements $P = P_1\biguplus \cdots \biguplus P_n$, $n\ge 2$, such that 
\begin{itemize}
\item The set $P_n$ is an order filter containing a minimal element of $P$.

\item The set $P_n$ is a chain.

\item The partition $P_1\biguplus \cdots \biguplus P_{n-1}$ is a chain peeling of $P\setminus P_n$.
\end{itemize}
\end{definition}
In other words, if the order filter generated by a minimal element of
$P$ is a chain, we may ``peel off'' this chain and 
obtain a smaller poset with one minimal element less. If this process can go on until the empty
poset is all that remains, $P$ has a chain peeling.

A poset that has a chain peeling is {\em chain
  peelable}. Observe that if $P$ is chain peelable, then so is every
ideal of $P$. 

For root poset ideals, as we shall see, chain peelability
turns out to be equivalent to supersolvability. One direction of this
correspondence is fairly straightforward to establish.

\begin{lemma}\label{le:peel}
If $\I$ is a root poset ideal, any chain peeling of $\I$ is a supersolving partition of $\I$.
\begin{proof}
By induction on the rank, it is sufficient
to show that if $F_\I(\alpha)$ is a chain for some $\alpha\in \Delta$,
then this chain does not contain an entire $2$-flat.

Suppose $F_\I(\alpha)$ is a chain. Assume $\gamma > \gamma' \in
F_\I(\alpha)$ and that $\gamma, \gamma'$ are the two smallest elements
in $ F_\I(\alpha) \cap \spn\{\gamma, \gamma'\}$. Lemma
\ref{le:chainroot} shows that $\gamma - \gamma'$ is a positive
multiple of some $\beta\in \Phi^+$. Since $\beta < \gamma$, $\beta\in
\I \cap \spn\{\gamma, \gamma'\}$. By minimality of $\gamma$ and $\gamma'$,
$\beta\not\in F_\I(\alpha)$. Thus, $F_\I(\alpha)$ does not contain the
$2$-flat $\I \cap \spn\{\gamma, \gamma'\}$.
\end{proof}
\end{lemma}

\section{Simply laced root systems}\label{se:simply_laced}
For the moment, let us focus on simply laced root systems, leaving
multiple laces to the next section. The arguments in this section are roughly
classification independent. More precisely, we rely on the classification of
finite root systems for two facts about the Dynkin diagram, namely
that it contains no circuit and that if there is a node with degree at
least three, then the degree is precisely three and this node is unique
in its connected component. Both facts are easily shown from first
principles without consulting the full classification machinery.

When $\Phi$ is of type $D_4$, the ideal which consists of all elements
strictly below the sum of all simple roots has special properties. We
shall see that it is in some sense the unique minimal non-supersolvable ideal
in the simply laced setting. Let us introduce a name for ideals that
contain (an isomorphic copy of) it.

\begin{definition}\label{de:star}
An ideal $\I \subseteq \Phi^+$ is a {\em star ideal} if there are distinct
simple roots $\alpha_1, \alpha_2, \alpha_3, \alpha_4 \in \Delta$ such
that $\gamma_4 = \alpha_1 + \alpha_2 + \alpha_3$,
$\gamma_3 = \alpha_1 + \alpha_2 + \alpha_4$ and $\gamma_1 = \alpha_2 +
\alpha_3 + \alpha_4$ all are roots that belong to $\I$. 
\end{definition}

Suppose $\alpha_i$ are as in Definition \ref{de:star}. Using that all
roots have connected support and that there are no circuits 
in the Dynkin diagram, it is easy to see that $\alpha_1$, $\alpha_3$,
$\alpha_4$ must all be neighbours of $\alpha_2$. In particular,
$\alpha_1$, $\alpha_3$, $\alpha_4$ are pairwise orthogonal and
$\alpha_2$ has no other neighbours.

\begin{lemma} \label{le:nostar}
If $\I$ is a star ideal, then $\I$ is not supersolvable.
\begin{proof}
Assume in order to obtain a contradiction that $\I$ is a
star ideal with supersolving partition $\I = \Pi_1\biguplus \cdots \biguplus
\Pi_n$ and that $n$ is minimal among all supersolving partitions of
star ideals. Let $\gamma_i$ and $\alpha_j$ be as in Definition
\ref{de:star}. By Lemma \ref{le:twocases} and Lemma
\ref{le:ideal}, $\I_n \setminus \Pi_n$ is an ideal either in $\Phi^+$
or in the subsystem $\Phi'^+$ defined in Lemma \ref{le:ideal}. By
minimality of $n$, this ideal is not a star ideal.

Suppose first that we are in case (a) of Lemma \ref{le:twocases}. Then,
$\Pi_n = F_\I(\alpha_i)$ for some $i \in \{1,2,3,4\}$. Thus, at least two
of the incomparable $\gamma_i$ belong to $\Pi_n$, contradicting that
$F_\I(\alpha_i)$ is a chain.

If instead situation (b) of Lemma \ref{le:twocases} holds, we may
without loss of generality assume that either $\Pi_n = G_\I (\alpha_1,
\beta, a, b)$ for some $\beta\neq \alpha_i$, or else $\Pi_n = G_\I(\alpha_1,
\alpha_2, a, b)$. In the former case, the flat generated by $\gamma_3$
and $\gamma_4$ is contained in $\Pi_n$. In the latter, $\Pi_n$ instead
contains the flat generated by $\alpha_2$ and $\gamma_1$. Hence, $\Pi_n$
contains an entire $2$-flat, providing the contradiction.  
\end{proof}
\end{lemma}

If the subgraph of the Dynkin diagram induced by some simple roots
$\alpha_1,\ldots, \alpha_k\in \Delta$ is a path, then
$\alpha_1+\cdots+\alpha_k\in \Phi^+$ since, in particular, this
subgraph is connected. It is worthwile to introduce some notation pertaining to
roots of this form.

\begin{definition}\label{de:path}
Let $\path \Phi$ denote the set of positive roots $\gamma\in \Phi^+$
satisfying that $\supp \gamma$ is a path in the Dynkin diagram of $\Delta$
and $\gamma_\alpha = 1$ for every $\alpha \in \supp \gamma$.
\end{definition}

\begin{lemma}\label{le:path}
If $\Phi$ is simply laced and $\I$ is not
a star ideal, then $\I \subseteq \path\Phi$.
\begin{proof}
If $\I \not \subseteq \path \Phi$, there is some $\gamma\in \I \setminus \path
\Phi$ satisfying $\gamma' \in \path\Phi$ for all $\gamma' < \gamma$. For some
$\alpha\in \Delta$ and $\gamma' < \gamma$, we have $\gamma =
\alpha + \gamma'$. By Lemma \ref{le:rank2}, $\alpha \not \in \supp \gamma'$ since
$\Phi$ is simply laced and $\gamma'\in \path\Phi$. Therefore, the
support of $\gamma$ does not form a path in the Dynkin diagram of
$\Phi$, so this support must contain the diagram of $D_4$ as a
subgraph. Thus, $\I$ is a star ideal. 
\end{proof}
\end{lemma}

Since $\Phi$ is finite, the Dynkin diagram contains no circuits. If $\alpha_1, \alpha_2\in
\Delta$ belong to the same irreducible component, the diagram thus
contains exactly one path with $\alpha_1$, $\alpha_2$ as
endpoints. Let $(\alpha_1\frown\alpha_2) \in \path \Phi$ denote the
sum of the simple roots that comprise this path.

\begin{lemma}\label{le:nostarchain}
If $\Phi$ is simply laced and $\I$ is not a star ideal, then $\I$ is chain peelable.
\begin{proof}
Without loss of generality, assume $\Delta \subseteq \I$ and that
$\Phi$ is irreducible.

Suppose $\I$ is not a star ideal. By Lemma \ref{le:path},
$\I\subseteq \path \Phi$. Consider any leaf $\alpha\in \Delta$ in the
Dynkin diagram. If $F_\I(\alpha)$ is not a chain, there are $\alpha',
\alpha''\in \Delta$ such that $p' = (\alpha
\frown \alpha')\in \I$ and $p'' = (\alpha \frown \alpha'')\in \I$ are
incomparable. In particular, $p',p''\ge \alpha_2$ for some
$\alpha_2\in \Delta$ which has degree at least $3$ in the Dynkin
diagram. By the classification of simply laced diagrams, $\alpha_2$ must be the unique
node with degree $3$. Since $\I$ is not a star ideal, we may denote the
neighbours of $\alpha_2$ by $\alpha_1, \alpha_3, \alpha_4\in \Delta$,
in an appropriate order, so that $(\alpha_1 \frown \alpha_3) \not \in
\I$. 

If we now let $\alpha$ be the leaf satisfying $\alpha_1\in \supp
(\alpha \frown \alpha_2)$, we cannot find $\alpha', \alpha''\in \I$ as
above, since $p'$ or $p''$ would have to be larger than $(\alpha_1\frown
\alpha_3)$. Thus, $F_\I(\alpha)$ is a chain. The lemma follows by
induction on the rank.
\end{proof}
\end{lemma}

We are now in a position to characterise supersolvable root ideal
arrangements in simply laced types. In particular, Theorem
\ref{th:main_peelable} is established for such root systems.

\begin{theorem}\label{th:simply_laced}
If $\Phi$ is simply laced, the following assertions are equivalent for
an ideal $\I\subseteq \Phi^+$ in the root poset:
\begin{itemize}
\item[(a)] $\I$ is supersolvable.
\item[(b)] $\I$ is chain peelable.
\item[(c)] $\I$ is not a star ideal.
\end{itemize}
\begin{proof}
Lemma \ref{le:nostar} shows (a) $\Rightarrow$ (c), whereas (c) $\Rightarrow$
(b) is Lemma \ref{le:nostarchain}. Finally, the implication (b) $\Rightarrow$
(a) is established by Lemma \ref{le:peel}.
\end{proof}
\end{theorem}

In particular, all ideals in type $A$ root posets are supersolvable
whereas every root poset of type $D$ or $E$ contains a unique minimal star ideal,
containment of which characterises non-supersolvability.

\section{Multiply laced root systems}\label{se:multiply_laced}

In order to complete the proof of Theorem \ref{th:main_peelable}, we now turn
to root systems that are not simply laced. Here we make use of the
classification; the irreducible types that we need to consider are
$B_n$, $C_n$, $F_4$ and $G_2$. The only case which requires a little
care is type $F_4$.

\begin{theorem}\label{th:BCG}
If $\Phi$ is of type $B_n$, $C_n$ or $G_2$, every root poset ideal
$\I\subseteq \Phi^+$ is both chain peelable and supersolvable.
\begin{proof}
Every chain peelable ideal is supersolvable by Lemma
\ref{le:peel}. Since chain peelability is inherited by subideals, it
suffices to show that the entire root poset $\Phi^+$ is chain
peelable. In types $B_n$ and
$C_n$, the statement follows by induction on $n$. Namely, if $n\ge 3$,
it is enough to note that $F_{\Phi^+}(\alpha)$ is a chain if $\alpha$
is the simple root corresponding to the leaf which is not incident to the double bond in the
Dynkin diagram. In the base case $B_2 = C_2$, as well as in type
$G_2$, chain peelability is readily verified.
\end{proof}
\end{theorem}

Now consider the type $F_4$ root poset with simple roots labelled as in
Figure \ref{fi:F4}. Let $\IH$ be the ideal generated by the three
roots $\eta_1 = \alpha_1+2\alpha_2+\alpha_3$,
$\eta_2 = \alpha_1+\alpha_2+\alpha_3+\alpha_4$ and
$\eta_3 = 2\alpha_2+\alpha_3+\alpha_4$. Thus, $\IH$ consists of the positive
roots of height at most $4$.

Concluding the proof of Theorem \ref{th:main_peelable}, we now verify
that containment of $\IH$ characterises both supersolvability and
chain peelability in type $F_4$.

\begin{theorem}\label{th:F4}
Suppose $\Phi$ is of type $F_4$. For an ideal $\I\subseteq\Phi^+$ in
the root poset, the following are equivalent:
\begin{itemize}
\item[(a)] $\I$ is supersolvable.
\item[(b)] $\I$ is chain peelable.
\item[(c)] $\I \not \supseteq \IH$.
\end{itemize}
\begin{proof}
Since all roots of height $3$ are covered
by two roots each, it is clear that $F_\I(\alpha)$ is not a chain
for any  simple root $\alpha$ if $\I\supseteq \IH$. Hence such an
ideal cannot be chain peelable. Conversely, it is easy to see that the ideals
$\Phi^+\setminus F_{\Phi^+}(\eta_i)$ are chain
peelable. If $\I \not \supseteq \IH$, then $\I$ is contained in at
least one of these three ideals. Thus, $\I$ is chain peelable if and
only if $\I\not \supseteq \IH$.

Now suppose $\I \supseteq \IH$. It remains to show that $\I$ is not
supersolvable. Case (a) of Lemma \ref{le:twocases} was ruled out
above, so if a supersolving partition $\I = \Pi_1\biguplus \Pi_2
\biguplus \Pi_3 \biguplus \Pi_4$ were to exist, the only possibility
would be to have $\Pi_4 = G_\I(\alpha,
\beta, a, b)$ for suitable chosen $\alpha, \beta, a, b$. Labelling the
simple roots as in the Dynkin diagram in Figure \ref{fi:F4}, there are four
possible choices to consider. For each case we now provide a $2$-flat entirely
contained in $G_\I(\alpha, \beta, a, b)$, thereby ruling out the
possibility of a supersolving partition:
\begin{itemize}
\item The first possibility is $G_\I(\alpha_1, \alpha_2, 1, 1)$ which contains the $2$-flat
  consisting of $2\alpha_2 + \alpha_3$ and
  $\alpha_2+\alpha_3+\alpha_4$.
\item Second, $G_\I(\alpha_2, \alpha_3, 1, 1)$ contains the $2$-flat comprised of
  $\alpha_1 + \alpha_2$ and $2\alpha_2 + \alpha_3$.
\item Third, $G_\I(\alpha_3, \alpha_4, 1, 1)$ contains the $2$-flat
  consisting of $2\alpha_2 + \alpha_3$ and  $\alpha_1+\alpha_2+\alpha_3$.
\item Finally, the set $G_\I(\alpha_2, \alpha_3, 2, 1)$ contains the $2$-flat
  which consists of $\alpha_1 + \alpha_2 + \alpha_3$,
  $\alpha_2+\alpha_3+\alpha_4$ and (if $\I$ contains it)
  $\alpha_1+2\alpha_2+2\alpha_3+\alpha_4$.
\end{itemize}

\end{proof}
\end{theorem}

\section{Line-closedness, quadraticity and the Koszul property}
In this final section, we shall apply the classification of
supersolvable root ideal arrangements that was obtained in the
previous sections in order to establish that these are the only root
ideal arrangements whose Orlik-Solomon algebras are Koszul. For an
account of how the Koszul property of Orlik-Solomon algebras is
relevant to the study of arrangements, we refer to \cite{yuzvinsky}. A
general introduction to Koszul algebras can be found in \cite{froberg}.

In an effort to better understand what makes Orlik-Solomon algebras
quadratic, Falk \cite{falk} introduced the concept of line-closedness:
\begin{definition} Let $\A$ be an arrangement.
\begin{itemize}
\item A subset $S\subseteq \A$ is {\em $2$-closed} if $\cl (\{a,b\})
  \subseteq S$ for all $a,b\in S$.
\item The arrangement $\A$ is {\em line-closed} if every $2$-closed subset
of $\A$ is a flat of $\A$.
\end{itemize}
\end{definition}
Falk showed that every arrangement with a quadratic
Orlik-Solomon algebra is line-closed.\footnote{A counterexample to the conjectured converse was found by
Yuzvinsky; cf.\ \cite{DY}.} In particular, being line-closed is a
necessary condition for having a Koszul Orlik-Solomon algebra. We now
observe that line-closed ideals cannot contain the minimal
non-supersolvable ideals that were identified in Sections
\ref{se:simply_laced} and \ref{se:multiply_laced}. Recall from those
sections the definition of star ideals and the ideal $\IH$ in the
type $F_4$ root poset, respectively.

\begin{proposition}\label{pr:line_closed}
Suppose $\Phi$ is a finite root system and $\I \subseteq \Phi^+$ is an
ideal in the root poset.
\begin{itemize}
\item[(a)] If $\Phi$ is simply laced and $\I$ is a
  star ideal, then $\I$ is not line-closed.
\item[(b)] If $\Phi$ is of type $F_4$ and $\I \supseteq \IH$, then $\I$ is not line-closed. 
\end{itemize}
\begin{proof}
For (a), suppose $I$ is a star ideal. Let $\gamma_i$ and $\alpha_j$ be as in Definition
\ref{de:star}. Observe that the four roots $\alpha_2$, $\gamma_1$,
$\gamma_3$ and $\gamma_4$ form a $2$-closed subset of $\I$. Since
(for example) $\alpha_1 = \frac{1}{2}(\gamma_3 + \gamma_4 - \gamma_1 -
\alpha_2)$, this subset is not a flat. Hence, $\I$ is not line-closed. 

Turning to (b), let $\alpha_i$ and $\eta_j$ be as defined prior to
Theorem \ref{th:F4}. Suppose $\IH \subseteq \I$ and define
\[
S = \{\eta_1, \eta_2, \eta_3, \alpha_3\} \cup (\I \cap \{\alpha_3 +
\eta_3, \eta_1 + \eta_2\}).
\]
Then, $S$ is a $2$-closed subset of $\I$. However, $\alpha_2 =
\frac{1}{3}(\eta_1 - \eta_2 + \eta_3 - \alpha_3) \in \I
\setminus S$ so that $S$ is not a flat and $\I$ is not line-closed. 
\end{proof}
\end{proposition}

\begin{corollary}
For a root ideal arrangement $\A_\I$, the following assertions are equivalent:
\begin{itemize}
\item[(a)] $\A_\I$ is supersolvable.
\item[(b)] $\A_\I$ is line-closed.
\item[(c)] $\OS{\A_\I}$ is Koszul.
\end{itemize}
\begin{proof}
Proposition \ref{pr:line_closed} and Theorems
\ref{th:simply_laced}, \ref{th:BCG} and \ref{th:F4} show (ii)
$\Rightarrow$ (i). As was previously mentioned, (i) $\Rightarrow$ (iii)
and (iii) $\Rightarrow$ (ii) follow from results of Shelton and
Yuzvinsky \cite{SY} and Falk \cite{falk}, respectively.
\end{proof}
\end{corollary}

In particular, Theorem \ref{th:main_koszul} is established.

\end{document}